\documentclass{gtart}
\def\ifplaintex{\expandafter\ifx\csname documentclass\endcsname\relax}

\def\gtp{{\mathsurround=0pt\it $\cal G\mskip-2mu$eometry \&\ 
$\cal T\!\!$opology $\cal P\!$ublications}}  % GT publications

\def\recd{{\small Received:\qua\receiveddate\ifx\reviseddate\relax
\else\qquad Revised:\qua\reviseddate\fi\par}} 

\def\lognumber#1{\def\thelognumber{#1}}
\def\volumenumber#1{\def\thevolumenumber{#1}}
\def\volumeyear#1{\def\thevolumeyear{#1}}
\def\papernumber#1{\def\thepapernumber{#1}}
\def\pagenumbers#1#2{\def\startpage{#1}\def\finishpage{#2}}
\def\published#1{\def\publishdate{#1}}

\def\received#1{\def\receiveddate{#1}}
\def\revised#1{\def\reviseddate{#1}}
\def\accepted#1{\def\accepteddate{#1}}
\def\asciititle#1{\def\theasciititle{#1}}

\def\asciiaddress#1{\def\theasciiaddress{#1}}

\long\def\asciiabstract#1{\long\def\theasciiabstract{#1}}

\let\\\par\let\thelognumber\relax\let\thevolumenumber\relax
\let\thepapernumber\relax\let\thevolumeyear\relax\let\startpage\relax
\let\finishpage\relax\let\publishdate\relax\let\receiveddate\relax
\let\reviseddate\relax\let\accepteddate\relax\let\theasciititle\relax
\let\theasciiauthors\relax\let\theasciiaddress\relax
\let\theasciiabstract\relax

\let\theasciiemail\relax

\ifplaintex
\font\logobig=cmssbx10 scaled 3836
\font\logomed=cmssbx10 scaled 2557
\else
\font\logobig=cmssbx10 scaled 4200
\font\logomed=cmssbx10 scaled 2800
\fi

\long\def\makeagttitle{   %%% start of definition of \makeagttitle
\count0=\startpage
\agt\hfill      %   Journal title (top left) 
\hbox to 45truept{\vbox to 0pt{\vglue -13truept{\logomed A\kern -.37em{\logobig 
T}\kern -.38em G}\vss}\hss}
\break
{\small Volume \thevolumenumber\ (\thevolumeyear)
\startpage--\finishpage\nl
Published: \publishdate}

\vglue .25truein

{\parskip=0pt\leftskip 0pt plus
1fil\def\\{\par\smallskip}{\Large\bf\thetitle}\par\medskip} \vglue
0.05truein

{\parskip=0pt\leftskip 0pt plus 1fil\def\\{\par}{\sc\theauthors}
\par\medskip}%
 
\vglue 0.03truein 

{\small\leftskip 25truept\rightskip 25truept{\bf Abstract}\stdspace\theabstract

{\bf AMS Classification}\stdspace\theprimaryclass
\ifx\thesecondaryclass\relax\else; \thesecondaryclass\fi\par
{\bf Keywords}\stdspace \thekeywords\par}\vglue 7truept

}   %%%% end of definition of \makeagttitle

\ifplaintex
\font\phead=cmsl9 scaled 950
\font\pnum=cmbx10 scaled 913
\font\pfoot=cmsl9 scaled 950
\headline{\vbox to 0pt{\vskip -4.5mm\line{\small\phead\ifnum
\count0=\startpage ISSN 1472-2739 (on-line) 1472-2747 (printed)
\hfill {\pnum\folio}\else\ifodd\count0\def\\{ }% 
\ifx\theshorttitle\relax\thetitle\else\theshorttitle\fi\hfill{\pnum\folio}
\else\def\\{ and }{\pnum\folio}\hfill\ifx\theshortauthors\relax\theauthors
\else\theshortauthors\fi\fi\fi}\vss}}
\footline{\vbox to 0pt{\vglue 0mm\line{\small\pfoot\ifnum\count0=\startpage
\copyright\ \gtp\hfill\else
\agt, Volume \thevolumenumber\ (\thevolumeyear)\hfill\fi}\vss}}
\else
\font=cmsl9 scaled 1050
\font\lnum=cmbx10 
\font=cmsl9 scaled 1050
\makeatletter
\def\@oddhead{{\small\lhead\ifnum\count0=\startpage ISSN 1472-2739 
(on-line) 1472-2747 (printed)\hfill {\lnum\number\count0}\else\ifodd\count0
\def\\{ }\ifx\theshorttitle\relax \thetitle \else\theshorttitle\fi\hfill
{\lnum\number\count0}\else\def\\{ and }{\lnum\number\count0}
\hfill\ifx\theshortauthors\relax 
\theauthors\else\theshortauthors\fi\fi\fi}}\def\@evenhead{\@oddhead}
\def\@oddfoot{\small\lfoot\ifnum\count0=\startpage\copyright\ \gtp\hfill\else
\agt, Volume \thevolumenumber\ (\thevolumeyear)\hfill\fi}
\def\@evenfoot{\@oddfoot}
\makeatother
\fi
\let\maketitlepage\makeagttitle
\let\makeshorttitle\maketitlepage
\let\maketitle\maketitlepage

\newwrite\gtoutfile
\long\gdef\makeheadfile{  %%% start of definition of \makeheadfile
{\def\\{, }\def\s{ }
\immediate\openout\gtoutfile head.xxx
\immediate\write\gtoutfile{To: math@arxiv.org}
\immediate\write\gtoutfile{Subject: put OR rep NNNNN:ppppp}
\immediate\write\gtoutfile{--text follows this line--}
\immediate\write\gtoutfile{Proxy-for: \ifx\theasciiauthors\relax
\theauthors\else\theasciiauthors\fi\s<\ifx\theasciiemail\relax\theemail\else\theasciiemail\fi>}
\immediate\write\gtoutfile{\noexpand\\}
\immediate\write\gtoutfile{Authors: \ifx\theasciiauthors\relax
\theauthors\else\theasciiauthors\fi}
{\def\\{ }\immediate\write\gtoutfile{Title: \ifx\theasciititle\relax
\thetitle\else\theasciititle\fi}}
\immediate\write\gtoutfile{Subj-class: GT or SG, GR etc}
\immediate\write\gtoutfile{MSC-class: \theprimaryclass\ifx\thesecondaryclass\relax\else, \thesecondaryclass\fi}
\immediate\write\gtoutfile{Journal-ref: Algebr. Geom. Topol. \thevolumenumber\s
(\thevolumeyear) \startpage-\finishpage}
\immediate\write\gtoutfile{Comments: Published by Algebraic and
Geometric Topology at}
\immediate\write\gtoutfile{\s\s\s  http://www.maths.warwick.ac.uk/agt/AGTVol\thevolumenumber/agt-\thevolumenumber-\thepapernumber.abs.html}
\immediate\write\gtoutfile{\noexpand\\}
\immediate\write\gtoutfile{}
\ifx\theasciiabstract\relax
\immediate\write\gtoutfile{\theabstract}\else
\immediate\write\gtoutfile{\theasciiabstract}\fi
\immediate\write\gtoutfile{}
\immediate\write\gtoutfile{\noexpand\\}
\immediate\write\gtoutfile{}
\immediate\closeout\gtoutfile}}  %%% end of definition of \makeheadfile

\def\maketitlepage{\makeagttitle\makeheadfile}
\let\makeshorttitle\maketitlepage
\let\maketitle\maketitlepage

\def\ifplaintex{\expandafter\ifx\csname documentclass\endcsname\relax}

\def\gtp{{\mathsurround=0pt\it $\cal G\mskip-2mu$eometry \&\ 
$\cal T\!\!$opology $\cal P\!$ublications}}  % GT publications

\def\recd{{\small Received:\qua\receiveddate\ifx\reviseddate\relax
\else\qquad Revised:\qua\reviseddate\fi\par}} 

\def\lognumber#1{\def\thelognumber{#1}}
\def\volumenumber#1{\def\thevolumenumber{#1}}
\def\volumeyear#1{\def\thevolumeyear{#1}}
\def\papernumber#1{\def\thepapernumber{#1}}
\def\pagenumbers#1#2{\def\startpage{#1}\def\finishpage{#2}}
\def\published#1{\def\publishdate{#1}}

\def\received#1{\def\receiveddate{#1}}
\def\revised#1{\def\reviseddate{#1}}
\def\accepted#1{\def\accepteddate{#1}}
\def\asciititle#1{\def\theasciititle{#1}}

\def\asciiaddress#1{\def\theasciiaddress{#1}}

\long\def\asciiabstract#1{\long\def\theasciiabstract{#1}}

\let\\\par\let\thelognumber\relax\let\thevolumenumber\relax
\let\thepapernumber\relax\let\thevolumeyear\relax\let\startpage\relax
\let\finishpage\relax\let\publishdate\relax\let\receiveddate\relax
\let\reviseddate\relax\let\accepteddate\relax\let\theasciititle\relax
\let\theasciiauthors\relax\let\theasciiaddress\relax
\let\theasciiabstract\relax

\let\theasciiemail\relax

\ifplaintex
\font\logobig=cmssbx10 scaled 3836
\font\logomed=cmssbx10 scaled 2557
\else
\font\logobig=cmssbx10 scaled 4200
\font\logomed=cmssbx10 scaled 2800
\fi

\long\def\makeagttitle{   %%% start of definition of \makeagttitle
\count0=\startpage
\agt\hfill      %   Journal title (top left) 
\hbox to 45truept{\vbox to 0pt{\vglue -13truept{\logomed A\kern -.37em{\logobig 
T}\kern -.38em G}\vss}\hss}
\break
{\small Volume \thevolumenumber\ (\thevolumeyear)
\startpage--\finishpage\nl
Published: \publishdate}

\vglue .25truein

{\parskip=0pt\leftskip 0pt plus
1fil\def\\{\par\smallskip}{\Large\bf\thetitle}\par\medskip} \vglue
0.05truein

{\parskip=0pt\leftskip 0pt plus 1fil\def\\{\par}{\sc\theauthors}
\par\medskip}%
 
\vglue 0.03truein 

{\small\leftskip 25truept\rightskip 25truept{\bf Abstract}\stdspace\theabstract

{\bf AMS Classification}\stdspace\theprimaryclass
\ifx\thesecondaryclass\relax\else; \thesecondaryclass\fi\par
{\bf Keywords}\stdspace \thekeywords\par}\vglue 7truept

}   %%%% end of definition of \makeagttitle

\ifplaintex
\font\phead=cmsl9 scaled 950
\font\pnum=cmbx10 scaled 913
\font\pfoot=cmsl9 scaled 950
\headline{\vbox to 0pt{\vskip -4.5mm\line{\small\phead\ifnum
\count0=\startpage ISSN 1472-2739 (on-line) 1472-2747 (printed)
\hfill {\pnum\folio}\else\ifodd\count0\def\\{ }% 
\ifx\theshorttitle\relax\thetitle\else\theshorttitle\fi\hfill{\pnum\folio}
\else\def\\{ and }{\pnum\folio}\hfill\ifx\theshortauthors\relax\theauthors
\else\theshortauthors\fi\fi\fi}\vss}}
\footline{\vbox to 0pt{\vglue 0mm\line{\small\pfoot\ifnum\count0=\startpage
\copyright\ \gtp\hfill\else
\agt, Volume \thevolumenumber\ (\thevolumeyear)\hfill\fi}\vss}}
\else
\font=cmsl9 scaled 1050
\font\lnum=cmbx10 
\font=cmsl9 scaled 1050
\makeatletter
\def\@oddhead{{\small\lhead\ifnum\count0=\startpage ISSN 1472-2739 
(on-line) 1472-2747 (printed)\hfill {\lnum\number\count0}\else\ifodd\count0
\def\\{ }\ifx\theshorttitle\relax \thetitle \else\theshorttitle\fi\hfill
{\lnum\number\count0}\else\def\\{ and }{\lnum\number\count0}
\hfill\ifx\theshortauthors\relax 
\theauthors\else\theshortauthors\fi\fi\fi}}\def\@evenhead{\@oddhead}
\def\@oddfoot{\small\lfoot\ifnum\count0=\startpage\copyright\ \gtp\hfill\else
\agt, Volume \thevolumenumber\ (\thevolumeyear)\hfill\fi}
\def\@evenfoot{\@oddfoot}
\makeatother
\fi
\let\maketitlepage\makeagttitle
\let\makeshorttitle\maketitlepage
\let\maketitle\maketitlepage

\newwrite\gtoutfile
\long\gdef\makeheadfile{  %%% start of definition of \makeheadfile
{\def\\{, }\def\s{ }
\immediate\openout\gtoutfile head.xxx
\immediate\write\gtoutfile{To: math@arxiv.org}
\immediate\write\gtoutfile{Subject: put OR rep NNNNN:ppppp}
\immediate\write\gtoutfile{--text follows this line--}
\immediate\write\gtoutfile{Proxy-for: \ifx\theasciiauthors\relax
\theauthors\else\theasciiauthors\fi\s<\ifx\theasciiemail\relax\theemail\else\theasciiemail\fi>}
\immediate\write\gtoutfile{\noexpand\\}
\immediate\write\gtoutfile{Authors: \ifx\theasciiauthors\relax
\theauthors\else\theasciiauthors\fi}
{\def\\{ }\immediate\write\gtoutfile{Title: \ifx\theasciititle\relax
\thetitle\else\theasciititle\fi}}
\immediate\write\gtoutfile{Subj-class: GT or SG, GR etc}
\immediate\write\gtoutfile{MSC-class: \theprimaryclass\ifx\thesecondaryclass\relax\else, \thesecondaryclass\fi}
\immediate\write\gtoutfile{Journal-ref: Algebr. Geom. Topol. \thevolumenumber\s
(\thevolumeyear) \startpage-\finishpage}
\immediate\write\gtoutfile{Comments: Published by Algebraic and
Geometric Topology at}
\immediate\write\gtoutfile{\s\s\s  http://www.maths.warwick.ac.uk/agt/AGTVol\thevolumenumber/agt-\thevolumenumber-\thepapernumber.abs.html}
\immediate\write\gtoutfile{\noexpand\\}
\immediate\write\gtoutfile{}
\ifx\theasciiabstract\relax
\immediate\write\gtoutfile{\theabstract}\else
\immediate\write\gtoutfile{\theasciiabstract}\fi
\immediate\write\gtoutfile{}
\immediate\write\gtoutfile{\noexpand\\}
\immediate\write\gtoutfile{}
\immediate\closeout\gtoutfile}}  %%% end of definition of \makeheadfile

\def\maketitlepage{\makeagttitle\makeheadfile}
\let\makeshorttitle\maketitlepage
\let\maketitle\maketitlepage

\lognumber{31}
\volumenumber{1}
\volumeyear{2001}
\papernumber{31}
\pagenumbers{605}{625}
\received{18 December 2000}
\revised{23 October 2001}
\accepted{24 October 2001}
\published{29 October 2001}

\usepackage{amsmath,amssymb,amscd}

\theoremstyle{plain}
\newtheorem{Thm}{Theorem}
\newtheorem{Cor}{Corollary}
\newtheorem{Prop}{Proposition}[section]
\newtheorem{Lem}[Prop]{Lemma}

\input{epsf}
\theoremstyle{definition}
\newtheorem{Def}{Definition}

\newtheorem*{remark}{Remark}

\begin{document}

\title{The Homflypt skein module of a connected\\sum of $3$-manifolds}
\asciititle{The Homflypt skein module of a connected\\sum of 3-manifolds}
\author{Patrick M. Gilmer\\Jianyuan K. Zhong}
\email{gilmer@math.lsu.edu, kzhong@coes.LaTech.edu}
\address{Department of Mathematics,   Louisiana State University\\Baton 
Rouge, LA 70803, USA}
\secondaddress{Program of Mathematics and Statistics\\Louisiana 
Tech University, Ruston, LA 71272, USA}
\asciiaddress{Department of Mathematics, Louisiana State University\\Baton 
Rouge, LA 70803, USA\\and\\Program of Mathematics and Statistics, Louisiana 
Tech University\\Ruston, LA 71272, USA}

\begin{abstract} If $M$ is an oriented $3$-manifold,  let $S(M)$ denote the Homflypt skein module of $M.$  We show that $S(M_1\# M_2)$ is isomorphic to  $S(M_1)\otimes S(M_2)$ modulo torsion. In fact, we show that
$S(M_1\# M_2)$ is isomorphic to $S(M_1)\otimes S(M_2)$ if we are working over a certain localized ring.  We show the similar result holds for relative skein modules.  If $M$ contains a separating $2$-sphere, we give conditions under which certain relative skein modules
of $M$ vanish over  specified localized rings. \end{abstract}

\asciiabstract{If M is an oriented 3-manifold, let S(M) denote
the Homflypt skein module of M. We show that S(M_1 connect sum M_2) is
isomorphic to S(M_1) tensor S(M_2) modulo torsion. In fact, we show
that S(M_1 connect sum M_2) is isomorphic to S(M_1) tensot S(M_2) if we are
working over a certain localized ring.  We show the similar result
holds for relative skein modules.  If M contains a separating
2-sphere, we give conditions under which certain relative skein
modules of M vanish over specified localized rings.}


\maketitlepage
\section{Introduction} 

We will be working with framed oriented links. By this we mean links equipped with a string orientation together with a nonzero normal vector field up to homotopy. The links described by figures in this paper will be assigned the ``blackboard'' framing which points to the right when travelling along an oriented strand.

\begin{Def}
{\bf The Homflypt skein module}\qua Let $k$ be a commutative ring containing $x^{\pm 1},$ $v^{\pm 1},$ $s^{\pm 1},$ and $\frac 1 {s-s^{-1}}.$ Let $M$ be an oriented $3$-manifold. The Homflypt skein module of $M$ over $k,$ denoted by $S_k(M)$, is the $k$-module freely generated by isotopy classes of framed oriented links in $M$ including the empty link, quotiented by the Homflypt skein relations given in the following figure.
$$
x^{-1}\quad\raisebox{-3mm}{\epsfxsize.3in\epsffile{left.ai}}\quad
        -\quad x\quad\raisebox{-3mm}{\epsfxsize.3in\epsffile{right.ai}}\quad
=\quad (\ s - \ s^{-1})\quad\raisebox{-3mm}{\epsfxsize.3in\epsffile{parra.ai}}\quad ,
$$
$$
\raisebox{-3mm}{\epsfxsize.35in\epsffile{framel.ai}}
        \quad=\quad (xv^{-1})\quad \raisebox{-3mm}{\epsfysize.3in\epsffile{orline.ai}}\quad  ,
$$
$$
L\ \sqcup \raisebox{-2mm}{\epsfysize.3in\epsffile{unknot.ai}}\quad
                                        =\quad {\frac{v^{-1}-v} {\ s - \ s^{-1}}}\quad L\quad.
$$
The last relation follows from the first two in the case $L$ is nonempty.\end{Def} 

\remark  (1) An embedding $f: M\to N$ of $3$-manifolds induces a well defined homomorphism $f_{*}: S_k(M)\to S_k(N)$. (2) If $N$ is obtained by adding a $3$-handle to $M$, the embedding $i: M\to N$ induces an isomorphism $i_{*}: S_k(M)\to S_k(N)$. (3) If $N$ is obtained by adding a $2$-handle to $M$, the embedding $i: M\to N$ induces an epimorphism $i_{*}: S_k(M)\to S_k(N)$. (4) If $M_1 \sqcup M_2$ is the disjoint union of $3$-manifolds $M_1$ and $M_2$, then $S_k(M_1 \sqcup M_2)\cong S_k(M_1)\otimes S_k(M_2)$.

Associated to a
partition of $n$, $\lambda=(\lambda_1\geq \dots \lambda_p \geq
1)$, $\lambda_{1}+\dots +\lambda_{p }=n$, is associated a Young
diagram of size $|\lambda|=n$, which we denote also by $\lambda$.
This diagram has $n$ cells indexed by \{$(i, j), \ 1\leq i\leq p$,
$1\leq j\leq \lambda_i$\}. If $c$ is the cell of index $(i, j)$ in
a Young diagram $\lambda$, its content
$cn(c)$ is defined by $cn(c)=j-i.$ Define
$$c_{\lambda,\mu}= 
v(s^{-1}-s)\sum_{c \in \mu}s^{-2cn(c)} +
v^{-1}(s-s^{-1})\sum_{c \in \lambda}s^{2cn(c)}$$

Let ${ \cal I}$ denote the submonoid  of the multiplicative monoid of 
${\mathbb{Z}}[v,s]$ generated by $v$, $s$,${s^{2n}-1}$ for all integers $n>0,$ and $c_{\lambda,\mu}$ for all pairs of Young diagrams $\lambda$, and $\mu,$ with $|\lambda| = |\mu|$, and $|\mu| \ne 0 .$ Let ${\cal {R}}$ be  ${\mathbb{Z}}[v,s]$ localized at  ${ \cal I}.$ \cite[7.2]{J}

\begin{Thm}
 $$S_{{\cal {R}}[x, x^{-1}]}(M_1 \# M_2)\cong S_{{\cal {R}}[x, x^{-1}]}(M_1)\otimes S_{{\cal {R}}[x, x^{-1}]}(M_2).$$
\end{Thm}
\noindent \remark J. Przytycki has proved the analog of this result for the Kauffman bracket skein module \cite{JP}. Our proof  follows the same general outline. We thank J. Przytycki for suggesting the problem of obtaining a similar result for the Homflypt skein module.

Let ${ \cal I}'$ denote the submonoid  of the multiplicative monoid of 
$\cal {R}$ generated by $v^4-s^{2n},$  for all $n.$ 
Let ${\cal {R}}'$ be  ${\cal {R}}$ localized at  ${ \cal I}'.$
It follows from \cite{ZG99}, $S_{{\cal {R}}'[x, x^{-1}]}(S^1\times S^2)$ is the free 
${\cal {R}}'[x, x^{-1}]$-module generated by the empty link. 

\begin{Cor}
 $S_{{\cal {R}}'[x, x^{-1}]}( \#^m S^1\times S^2)$  is a free module generated by the empty link.
\end{Cor}
\noindent \remark This allows us to define a ``Homflypt rational function'' $f$ in ${\cal {R}}'$ for oriented framed links in $\#^m S^1\times S^2$. If $L$ is such a link, one defines $f(L)$ by $L=f(L)\phi\in S_{{\cal {R}}'}( \#^m S^1\times S^2)$. A specific example is given in section 5.

Let  $l= \cal R$ with $x=v,$ then $S_l(M)$ is a version of the Homflypt skein module for unframed links. The next  two corollaries follows  from the universal coefficient property for skein modules which has been described by J. Przytycki \cite{JP}  for the Kauffman bracket skein module. The proof given there holds generally for essentially any skein module.

\begin{Cor}
  $S_l(M_1 \# M_2)\cong S_l(M_1)\otimes S_l(M_2).$\end{Cor}

Let  $l'= {\cal R}'$ with $x=v.$

\begin{Cor}
 $S_{l'}( \#^m S^1\times S^2)$  is a free $l'$-module generated by the empty link.
\end{Cor}

\begin{Def}
{\bf The relative Homflypt skein module}\qua Let $X=\{x_1, x_2,\cdots$, $x_{n}\}$ be a finite set of input framed points in $\partial M$, and let $Y=\{y_1, y_2,\cdots, y_{n}\}$ be a finite set of output framed points in the boundary $\partial M$. Define the relative skein module $S_k(M, X, Y)$ to be the $k$-module generated by relative framed oriented links in $(M,\partial M)$ such that $L \cap \partial M=\partial L=\{x_i, y_i\}$ with the induced framing, considered up to an ambient isotopy fixing $\partial M$, quotiented by the Homflypt skein relations. \end{Def}

Let $S(M)$ denote $S_{{\cal R}[x, x^{-1}]}(M,)$ and let $S(M,X, Y)$ denote $S_{{\cal R}[x, x^{-1}]}(M,X, Y).$  
We have the following version of Theorem 1 for relative skein modules.
At this point we must work over the field of fractions of ${\mathbb{Z}}[x,v,s]$ which we denote by ${\cal {F}}.$   
This is because we do not know whether the relative skein module of a handlebody is free. We conjecture that it is free. In the proof of Theorem 1, we use the absolute case first obtained by  Przytycki \cite{JP90}. We state Theorem 2 over ${\cal F},$ but conjecture it over ${\cal {R}}'[x, x^{-1}].$

\begin{Thm}
Let $M_1$ and $M_2$ be connected oriented 3-manifolds.
Let  $X_i=\{x_1, x_2,\cdots, x_{n}\}$ be a finite set of input framed points in $\partial M_i$, and let $Y_i=\{y_1, y_2,\cdots, y_{n}\}$ be a finite set of output framed points in the boundary $\partial M_i.$ Let $X = X_1 \cup X_2,$ and $Y = Y_1 \cup Y_2$, then
$$S_{\cal {F}}(M_1 \# M_2, X,Y)\cong S_{\cal {F}}(M_1, X_1,Y_1)\otimes S_{\cal {F}}(M_2,X_2,Y_2).$$
\end{Thm}

We also have the following related result. Let ${ \cal I}_r$ denote the submonoid  of the multiplicative monoid of 
${\mathbb{Z}}[x,v,s]$ generated by $x,$ $v,$  $s,$ ${s^{2n}-1}$ for all integers $n>0,$ and $x^r -1-c_{\lambda,\mu}$ for all pairs of Young diagrams $\lambda$, and $\mu,$ such that $|\lambda| - |\mu|=r$, and $|\mu|$ and $|\lambda|$ are not both zero.  Let $k_r$ be  ${\mathbb{Z}}[x,v,s]$ localized at ${ \cal I}_r.$ Note $k_0= {\cal R}[x, x^{-1}].$

\begin{Thm} Suppose $M$ is connected and contains a  2-sphere $\Sigma,$
such that $M - \Sigma$ has two connected components. Let $M'$ be one of these components. If  
$|X \cap M'|-|Y \cap M'|=r \ne 0,$ then  $S_{k_r}(M,X, Y)=0.$ 
\end{Thm}

In section 2, we prove that there is an epimorphism from $S({\cal H}_{m_1})\otimes S({\cal H}_{m_2})$
to $S({\cal H}_{m_1} \# {\cal H}_{m_2}).$ Here and below, we let ${\cal H}_m$
denote a handlebody of genus $m.$ 
In section 3, we prove Theorem 1 in the case of handlebodies. We prove Theorem 1 in the general case in section 4. Section 5 describes the class of a certain link  
in the $ S^1\times S^2 \# S^1\times S^2.$ Section 6 gives a proof of a lemma needed in section 2. In section 7, we discuss the proofs of Theorems 2 and 3.  

\section{Epimorphism for Handlebodies}

\subsection{The $n$th Hecke algebra of Type A}
We will use the related work  of  C. Blanchet \cite{cB98}, A. Aiston and H. Morton  \cite{AM98} on the $n$th Hecke algebra of Type A.  This is summarized in section $3$ of \cite{ZG99} whose conventions we follow. For the convenience of the reader, we give the basic definitions in this subsection. 

Note that  $s^{2n}-1$ is invertible in ${\cal R}$ for integers $n>0.$ It follows that the quantum integers $[n]=\dfrac{s^n-s^{-n}}{s-s^{-1}}$ for $n> 0$ are invertible in $k$. Let $[n]!=\prod_{j=1}^{n}[j]$, so $[n]!$ is invertible for $n>0$.

\medskip
{\bf The Hecke category}\qua The $k$-linear Hecke category $H$ is defined as follows. An object in this category is a disc $D^2$ equipped with a set of framed points. If $\alpha=(D^2, l)$ and $\beta=(D^2, l')$ are two objects, the module $Hom_H(\alpha, \beta)$ is ${\cal{S}}(D^2 \times [0,1], l\times 1, l'\times 0)$. The notation $H(\alpha, \beta)$ and $H_{\alpha}$ will be used for $Hom_H(\alpha, \beta)$  and $H(\alpha, \alpha)$ respectively. The composition of morphisms are by stacking the first one on the top of the second one.

Let $\otimes$ denote the monoid structure on $H$ given by embedding two disks $D^2$ side by side into one disk. For a Young diagram $\lambda$, by assigning each cell of ${\lambda}$ a point equipped with the horizontal (to the left) framing, we obtain an object of the category $H$ denoted by $\square_{\lambda}$. When $\lambda$ is the Young diagram with a single row of $n$ cells, $H_{\square_{\lambda}}$ will be denoted by $H_n$, which is the $n$th Hecke algebra of type A \cite{MT90}, \cite{vT90}.

For each permutation $\pi\in S_n$,  a positive permutation braid, $w_{\pi}$, is a braid which realizes
the permutation $\pi$ with all crossings positive \cite{M93}.
Let $\sigma_i \in H_n,\ i=1,\dotsc,n-1$, be the positive permutation corresponding to the transposition $(i\ i+1)$. As in  \cite{AM98}, define
$$f_n=\dfrac{1}{[n]!}s^{-\frac{n(n-1)}{2}}\sum_{\pi \in S_n}(xs^{-1})^{-l(\pi)}\omega_{\pi}$$
and
$$g_n=\dfrac{1}{[n]!}s^{\frac{n(n-1)}{2}}\sum_{\pi \in S_n}(-xs)^{-l(\pi)}\omega_{\pi}$$
Here $l(\pi)$ is the length of $\pi$.  

\medskip
{\bf Idempotents in the Hecke Algebra} \cite{AM98}\qua
 For a Young diagram $\lambda$ of size $n$, let $F_{\lambda}$ be the element in $H_{\square_{\lambda}}$ formed with one copy of $[\lambda_i]!f_{\lambda_i}$ along the row $i$, for $i=1,\dotsc,p.$  We let $\lambda^{\vee}$ denote the Young diagram whose rows are the columns of $\lambda.$  Let $G_{\lambda}$ be the element in $H_{\square_{\lambda}}$ formed with one copy of $[\lambda_j^{\vee}]!g_{\lambda_j^{\vee}}$ along the column $j$, for $j=1,\dotsc,q$. 
Let $\Tilde{y}_{\lambda}=F_{\lambda}G_{\lambda}$, then $\Tilde{y}_{\lambda}$ is a quasi-idempotent. Let ${y}_{\lambda}$ be the normalized idempotent from $\Tilde{y}_{\lambda}$.

\medskip
{\bf A Basis for the ${n}$th Hecke Algebra ${H_n}$}\qua 
A standard tableau $t$ with the shape of a Young diagram $\lambda=\lambda(t)$ is a labeling of the cells, with the integers $1$ to $n$ increasing along the rows and the columns. Let $t'$ be the tableau obtained by deleting the cell numbered by $n$. Note the cell numbered by $n$ in a standard tableau is an {\em extreme cell}. C. Blanchet defines $\alpha_t \in H(n, \square_{\lambda})$ and $\beta_t \in H(\square_{\lambda}, n)$ inductively by
$$
\begin{matrix}
\alpha_1=\beta_1=1_1,\\
\alpha_t=(\alpha_{t'}\otimes 1_1)\rho_{t}y_{\lambda},\\
\beta_t=y_{\lambda}\rho_{t}^{-1}(\beta_{t'}\otimes 1_1).
\end{matrix}
$$
Here $\rho_{t} \in H(\square_{\lambda(t')}\otimes 1,\square_{\lambda})$ is the isomorphism given by an arc joining the added point to its place in $\lambda$ in the standard way. 

Note that $\beta_{\tau}\alpha_t=0$ if $\tau \ne t$, and $\beta_t \alpha_t=y_{\lambda(t)}$.

\begin{Thm}[Blanchet] The family $\alpha_t\beta_{\tau}$ for all standard tableaux $t, \tau$ such that $\lambda(t)=\lambda(\tau)$ for all Young diagrams $\lambda$ with $|\lambda|=n$ forms a basis for $H_n$.
\end{Thm}

Let ${\overline{H_n}}$ denote $H_n$ with the reversed string orientation.
$$\raisebox{-10mm}{\epsfxsize.0in
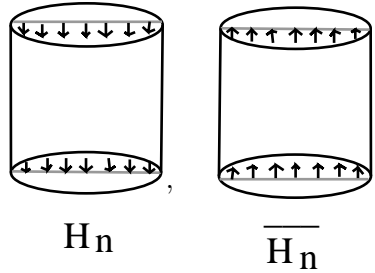
\epsffile{hnxy.eps}},\quad \ \raisebox{-12mm}{\epsfxsize.0in\epsffile{hnxy1.eps}}$$

\subsection{The Epimorphism on the Handlebodies}

If $X_m$ is a set of $m$ distinguished framed points in $D^2\times \{1\}$ and $Y_m$ be a set of $m$ distinguished  framed points in $D^2\times \{0\}$,
Let $\begin{matrix} = \\{(m)} \end{matrix}$ denote equality in $S(D^2\times I,X,Y)$ modulo the submodule ${\cal L} (m)$ generated by links which intersects $D^2\times \{\dfrac{1}{2}\}$ in less than  $m$ points.

In section 6, we derive:

\begin{Lem} Let $\lambda$, $\mu$ be two Young diagrams, and $m=|\lambda|+|\mu|.$
$$x^{2(|\mu|-|\lambda|)}\quad \raisebox{-18mm}{\epsfxsize.0in\epsffile{albesi01.eps}}\quad - \quad
\raisebox{-18mm}{\epsfxsize.0in\epsffile{out.eps}}\quad
\begin{matrix} = \\{(m)} \end{matrix} \quad c_{\lambda,\mu} \quad
 \raisebox{-10mm}{\epsfxsize.0in\epsffile{albesi02.eps}}$$
\end{Lem}

Let ${\cal H }_m$ be a handlebody of genus $m$. Let $D$ be a separating meridian disc of ${\cal{H}}_m$, let $\gamma=\partial{D}$. Let $({\cal{H}}_m)_{\gamma}$ be the manifold obtained by adding a $2$-handle to ${\cal{H}}_{m}$  along $\gamma$.
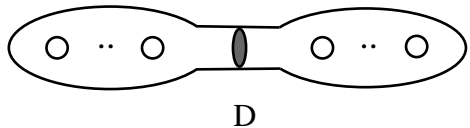
$$\raisebox{-3mm}{\epsfxsize.0in\epsffile{hm.eps}}$$
Let $V_D=[-1, 1]\times D$ be the regular neighborhood of $D$ in ${\cal{H}}_m$, $V_D$ can be projected into a disc $D_p=[-1, 1]\times [0, 1]$.

\begin{Lem}[The Epimorphism Lemma] The embedding $i:\ {\cal{H}}_m-D \to ({\cal{H}}_m)_{\gamma}$ induces an epimorphism: 
$$i_{*}: S({\cal{H}}_m-D)\twoheadrightarrow S(({\cal{H}}_m)_{\gamma}).$$
\end{Lem}

\begin{proof} Let $z_n$ be a link in ${\cal{H}}_m$ in general position with $D$ and cutting $D\ 2n$ times, let $z_n'=z_n\cap V_D$, i.e.,
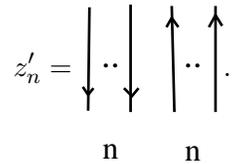
$$z_n'=\raisebox{-12.5mm}{\epsfxsize.0in\epsffile{zk.eps}}.$$
Note $z_n' \in H_n\otimes \overline{H_n}$. Using the basis elements $\alpha_t\beta_{\tau}$ of $H_n$ given in the previous theorem, $z_n'$ can be written as a linear combination of the elements $\alpha_t\beta_{\tau}\otimes {\alpha_{\sigma}\beta_{s}}^{*}$, where ${\alpha_{\sigma}\beta_{s}}^{*}$ is $\alpha_{\sigma}\beta_{s}$ with the reversed orientation. A diagram of $\alpha_t\beta_{\tau}\otimes {\alpha_{\sigma}\beta_{s}}^{*}$ is given by the following:
$$\alpha_t\beta_{\tau}\otimes {\alpha_{\sigma}\beta_{s}}^{*}=\raisebox{-18mm}{\epsfxsize.0in
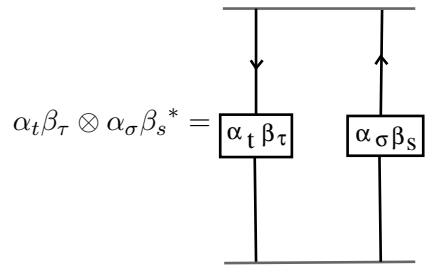
\epsffile{albesi.eps}}$$
By the inductive definition of $\alpha_t,\ \beta_{\tau},\ \alpha_{\sigma}, \ \beta_{s}$, an alternative diagram of $\alpha_t\beta_{\tau}\otimes {\alpha_{\sigma}\beta_{s}}^{*}$ is given by:
$$\alpha_t\beta_{\tau}\otimes {\alpha_{\sigma}\beta_{s}}^{*}=\raisebox{-18mm}{\epsfxsize.0in
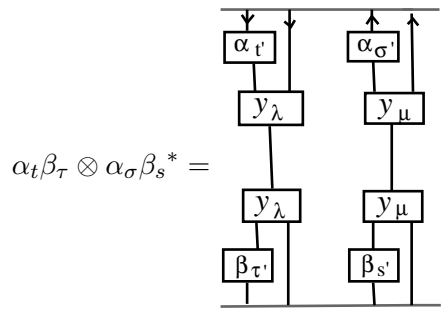
\epsffile{albesia.eps}}$$

We will consider the sliding relation given by: 
\begin{equation}
\raisebox{-12.5mm}{\epsfxsize.0in\epsffile{zno.eps}}\equiv \raisebox{-21.5mm}{\epsfxsize.0in\epsffile{zno1.eps}}\tag{I}
\end{equation}\vglue -10mm
From the above observation, we will be interested in the following relation:
\vglue -15mm\begin{equation}
\raisebox{-17mm}{\epsfxsize.0in\epsffile{albesi1.eps}}\quad \equiv \quad \raisebox{-17mm}{\epsfxsize.0in\epsffile{albesi2.eps}}\tag{II}
\end{equation}\vglue 2mm
From  Relation II, and Lemma 2.1, as $|\lambda|=|\mu|$, in
$S\left( (D^2 \times I)_{\gamma} \right)$ we have
$$c_{\lambda,\mu} \left( 
\alpha_t\beta_{\tau}\otimes { \alpha_{\sigma}\beta_{s} }^{*}
\right) \in 
{\cal L}_{|\lambda|+|\mu|}.$$ As $c_{\lambda,\mu}$ is invertible in $\cal R$,
we have that $\alpha_t\beta_{\tau}\otimes { \alpha_{\sigma}\beta_{s} }^{*}
\in {\cal L}(|\lambda|+|\mu|)$.
By induction, we can eliminate all elements of $({\cal{H}}_m)_{\gamma}$ which cut the $2$-disk $D_{\gamma}$ non-trivially. Thus $i_{*}$ is an epimorphism.
\end{proof}

\section {Isomorphism for handlebodies}
Recall that $({\cal{H}}_m)_{\gamma}$ is obtained by adding a $2$-handle to ${\cal{H}}_{m}$ along $\gamma$. From \cite{ZG99} section 2, we have $S(({\cal{H}}_m)_{\gamma})\cong  S({\cal{H}}_{m})/ R$, where $R$ is the submodule of $S({\cal{H}}_m)$ given by the collection $\{\Phi '( z ) - \Phi ''( z ) \mid z \in S({\cal{H}}_m, A, B)\}$. Here $A, \ B$ are two points on $\gamma$, which decompose $\gamma$ into two intervals $\gamma '$ and $\gamma ''$, $z$ is any element of the relative skein module $S({\cal{H}}_{m}, A, B)$ with $A$ an input point and $B$ an output point, and $\Phi '( z )$ and $\Phi ''( z )$ are given by capping off with $\gamma '$ and $\gamma ''$, respectively, and pushing the resulting links back into ${\cal{H}}_{m}$.

Let $I_0$ be the submodule of $S({\cal{H}}_m)$ given by the collection $\{p_{D} ( L ) - L\sqcup O \mid L \in S({\cal{H}}_m)\}$, where $O$ is the unknot. Locally, we have the following diagram description.
$$p_{D}( L )=\raisebox{-11.5mm}{\epsfxsize.0in\epsffile{lo1.eps}}, \quad L\sqcup O=\raisebox{-11.5mm}{\epsfxsize.0in\epsffile{lo.eps}}$$

\begin{Lem}
 $R=I_0$.
\end{Lem}

\begin{proof}
First note $R \supseteq I_0$. We need only  show that  $R \subseteq I_0$. Let $\pi$ be the projection map $\pi: S({\cal{H}}_m)\to S({\cal{H}}_m)/I_0$.  We will show that $\pi(R)=0$ in $S({\cal{H}}_m)/I_{0}$, i.e. $R\subseteq I_0$. We show this by proving now that $\pi(\Phi '( z )) = \pi(\Phi ''( z ))$ for any $z\in S({\cal{H}}_m, A, B)$.

Recall that $V_D=[-1, 1]\times D$ is the regular neighborhood of $D$ in ${\cal{H}}_m$. Let $D_1=\{-1\}\times D$ and $D_2=\{1\}\times D$. Let $\gamma_{1}=\partial D_1$ and $\gamma_{2}=\partial D_1$, note $\gamma_{1}$ and $\gamma_{2}$ are parallel to $\gamma$. 
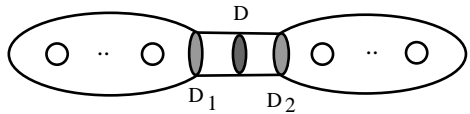
$$\raisebox{-12.5mm}{\epsfxsize.0in\epsffile{hmdd.eps}}$$
Let $I_{1}=\{p_{D_1} ( z ) - z\sqcup O \mid z \in S({\cal{H}}_m, A, B)\}$, $I_{2}=\{p_{D_2} ( z ) - z\sqcup O \mid z \in S({\cal{H}}_m, A, B)\}$, where locally
$$p_{D_1}( z )=\raisebox{-11.5mm}{\epsfxsize.0in\epsffile{lo11.eps}}, \quad p_{D_2}( z )=\raisebox{-11.5mm}{\epsfxsize.0in\epsffile{lo12.eps}},$$
$$ z\sqcup O=\raisebox{-11.5mm}{\epsfxsize.0in\epsffile{lo2.eps}}.$$
Let $\pi_{A,B}$ be the projection map $\pi_{A,B}: S({\cal{H}}_m, A, B) \to S({\cal{H}}_m, A, B)/(I_{1}+I_{2})$.  Note that $\Phi '( I_{i})=I_{0}$ and $\Phi ''( I_{i})=I_{0}$ for $i=1, 2$.

Let $z\in S({\cal{H}}_m-(D_1\cup D_2), A, B)$, then $\Phi '( z ) =\Phi ''( z )$ in $S({\cal{H}}_m-D_1-D_2)$, since $V_D=[-1, 1]\times D$ is a $3$-disc and closing a relative link along $\gamma$ by $\gamma'$ and $\gamma''$ in $V_D$ gives isotopic links. Let $\Phi:S({\cal{H}}_m-(D_1\cup D_2), A, B) \to
S({\cal{H}}_m-D_1-D_2)$ denote the map which sends $z$ to $\Phi '( z )= \Phi ''( z ).$

In general, let $z \in S({\cal{H}}_m, A, B)$. Now consider the following commutative diagram,
$$
\begin{CD}
S({\cal{H}}_m-(D_1\cup D_2), A, B)@>j_1>> 
S({\cal{H}}_m, A, B) @>\pi_{A,B}>>  S({\cal{H}}_m, A, B)/(I_{1}+I_{2})
\\
@V\Phi  VV   @V\Phi'  VV  @V{\bar \Phi'}  VV
\\
S({\cal{H}}_m-(D_1\cup D_2))@>j_2>> 
S({\cal{H}}_m) @>\pi>>  S({\cal{H}}_m)/I_{0}
\\
@A\Phi  AA   @A\Phi'' AA  @A{\bar \Phi''}  AA
\\
S({\cal{H}}_m-(D_1\cup D_2), A, B)@>j_1>> 
S({\cal{H}}_m, A, B) @>\pi_{A,B}>>  S({\cal{H}}_m, A, B)/(I_{1}+I_{2})
\end{CD}
$$

\noindent Here $j_1$ and $j_2$ are induced by inclusion maps. Also
$\bar \Phi',$ and $\bar \Phi''$ are induced by $\Phi',$ and $\Phi''$
respectively.
By an  argument similar to the proof of Lemma 2.2, the composition map $\pi_{A,B}j_1: \ S({\cal{H}}_m-(D_1\cup D_2), A, B) \to S({\cal{H}}_m, A, B)/(I_{1}+I_{2})$ is an epimorphism.

Take $z \in S({\cal{H}}_m, A, B)$, then $\pi_{A,B}(z)\in S({\cal{H}}_m, A, B)/(I_{1}+I_{2})$. As $\pi_{A,B}j_1$ is an epimorphism, there exists $z' \in S({\cal{H}}_m-(D_1\cup D_2), A, B)$ such that $\pi_{A,B}j_1(z')=\pi_{A,B}(z)$.  By the commutativity of the diagram, $\pi j_2(\Phi ( z' ))=\pi(\Phi '( z ))$ and 
$\pi j_2(\Phi ( z' ))=\pi(\Phi ''( z )).$ Thus $\pi(\Phi '( z )) = \pi(\Phi ''( z ))$.
\end{proof}

\begin{Cor}
The embedding ${\cal{H}}_{m} \to ({\cal{H}}_m)_{\gamma}$ induces an isomorphism
$$S({\cal{H}}_{m})/ I_0\cong S(({\cal{H}}_m)_{\gamma}).$$
\end{Cor}

Now we want to show that the embedding ${\cal{H}}_{m}-D \to ({\cal{H}}_m)_{\gamma}$ induces an isomorphism
$$S({\cal{H}}_{m}-D)\cong S(({\cal{H}}_m)_{\gamma}).$$

\begin{Lem}
$$S({\cal{H}}_m-D) \cap I_0=0.$$
\end{Lem}

\begin{proof}

Przytycki \cite{JP90} calculated the unframed Homflypt skein module of a handlebody. It follows from this, the universal coefficient property of skein modules and an  argument  of Morton in \cite{M93} section (6.2) that $S(H_m)$ is free.   As $S({\cal{H}}_{m}-D)$ is free,  the map  $S({\cal{H}}_{m}-D) \rightarrow S_{\cal F}({\cal{H}}_{m}-D),$ induced by 
${\cal R}[x,x^{-1}]\rightarrow {\cal F}$ is injective.  Let ${\cal {I}}_0=\{p_{D} ( L ) - L\sqcup O \mid L \in S_ {\cal {F}}({\cal{H}}_m)\}$. It is enough to show
 $S_{\cal F}({\cal{H}}_m-D) \cap {\cal I}_0=0.$

Let $\psi$ be the map from $S_{\cal F}({\cal{H}}_m) \rightarrow S_{\cal F}({\cal{H}}_m)$ given by $\psi(L)=p_D(L)-L\sqcup O$ for $L\in S_{\cal F}({\cal{H}}_m)$.  
$Image(\psi)={\cal I}_0$.

It also follows from
Przytycki's basis that the map induced by inclusion 
$S({\cal{H}}_{m}-D) \rightarrow S(H_m)$ is injective.
Let ${\cal B}_0$ be the image of a free basis for the module $S({\cal{H}}_{m}-D)$ 
in $S(H_m).$ ${\cal B}_0$ also a basis for injective image of
$S_{\cal F}({\cal{H}}_{m}-D)$ 
in $S_{\cal F}(H_m).$
Let $B_n$ be the subspace of $S_{\cal F}({\cal{H}}_m)$ generated by framed oriented links in ${\cal{H}}_m$ which intersect the disk $D \leq 2n$ times. Then we have a chain of vector spaces: 
$$B_0\subset B_1 \subset B_2 \subset \cdots \subset B_n \subset \cdots$$
${\cal{B}}_0$ is a basis for $B_0.$  The vector space $B_n/B_{n-1}$ is generated by elements of the form  $ {\alpha_t\beta_{\tau}\otimes {\alpha_{\sigma}\beta_{s}}^{*}}$
in a neighborhood of $D,$ where $|\lambda |= |\mu |$ is $n$. Let  ${\cal{B}}_n$ be a basis $B_n/B_{n-1},$ constructed by taking a maximal linearly independent subset  of the above generating set.  By the proof of Lemma 2.2, each  element of  ${\cal{B}}_n,$ where $n>0,$
 is an eigenvector for $\psi$ with nonzero  eigenvalue.
${\cal{B}}= \cup_{n\ge 0} {\cal{B}}_n$ is a basis for $S_{\cal F}({\cal{H}}_{m}).$
Let ${\cal{B}}'={\cal{B}}- {\cal{B}}_0$. Note $\psi({\cal B}_0)= 0$. So
${\cal I}_0= Image(\psi)= \psi({\cal B}').$

It follows that $\psi$ induces a one to one map: $B_n/B_{n-1}\rightarrow B_n/B_{n-1}$. Thus $\psi (<{\cal{B}}'>)\cap S_{\cal F}({\cal{H}}_m-D)=0$. The result follows.
\end{proof}

\begin{Thm}
The embedding ${\cal{H}}_{m}-D \to ({\cal{H}}_m)_{\gamma}$ induces an isomorphism
$$ S({\cal{H}}_m-D) \cong 
S(({\cal{H}}_m)_{\gamma})   .$$
\end{Thm}
\begin{proof}
From the above, we have the following commutative diagram
\[
\begin{CD}
0 @>>>  S({\cal{H}}_{m}-D)\cap I_0 @>>> S({\cal{H}}_{m}-D) @>>> S(({\cal{H}}_m)_{\gamma}) @>>> 0
\\
 @.   @VVV  @VVV @VVV
\\
0 @>>>  I_0 @>>> S({\cal{H}}_{m}) @>>> S(({\cal{H}}_m)_{\gamma}) @>>> 0
\end{CD}
\]
\end{proof}

${\cal{H}}_{m_1}\#{\cal{H}}_{m_2}$ is equal to ${\cal{H}}_{m_1 +m_2}$ with a $2$-handle added along the boundary of the meridian disc $D$ separating ${\cal{H}}_{m_1}$ from ${\cal{H}}_{m_2}$. Let $\gamma=\partial{D}$. Therefore we can consider ${\cal H}_1\# {\cal H}_2 =
({{\cal H}}_m)_{\gamma}.$ 
As ${\cal H}_{m_1 +m_2}- D= {\cal H}_{m_1}\sqcup {\cal H}_{m_2},$  the above theorem  says:

\begin{Cor} Let $B_1$ and $B_2$ denote the 3-balls we remove from ${\cal H}_{m_1}$ and ${\cal H}_{m_2}$ while forming ${\cal H}_{m_1}\# {\cal H}_{m_2}.$ The embedding  
$({\cal H}_{m_1} - B_1)\sqcup ({\cal H}_{m_2} - B_2) \rightarrow {\cal H}_{m_1}\# {\cal H}_{m_2}$ induces $$S({\cal H}_{m_1}) \otimes  S({\cal H}_{m_2}) \cong
S({\cal H}_{m_1}\# {\cal H}_{m_2})     $$
\end{Cor}

\section {The general case for absolute skein modules}

A connected oriented $3$-manifold with nonempty
boundary may be obtained from the handlebody ${\cal H}$ by adding some $2$-handles. If $M$ is closed, we will also need one 3-handle. As removing 3-balls from the interior of a 3-manifold does not change its Homflypt skein module, we may reduce Theorem 1 to the case that $M_1$ and $M_2$ are connected 3-manifolds with boundary. 
  
In this case, each $M_i$ is obtained from the handlebody ${\cal H}_{m_i}$ by adding some $2$-handles. 
Let $m=m_1+m_2$. Let $N$ be the manifold obtained by adding both sets of  $2$-handles to the boundary connected sum of  ${\cal H}_{m_1}$ and ${\cal H}_{m_2}$ which we identify with ${\cal H}_{m}$. Let $D$ be the disc in 
${\cal H}_{m}$ separating ${\cal H}_{m_1}$ from ${\cal H}_{m_2}.$  Let $\gamma=\partial{D},$ so ${\cal H}_{m_1}\# {\cal H}_{m_2}=({\cal H}_{m})_{\gamma}.$ Here and below  $P_{\delta}$ denotes the result of adding a 2-handle to a 3-manifold $P$ along a curve  $\delta$ in $\partial N.$
We can consider $M_1\# M_2$ as obtained from $({{\cal H}}_{m})_{\gamma}$ by adding those $2$-handles. Thus $N - D = M_1 \sqcup M_2,$  and $M_1\# M_2=N_{\gamma}.$
 
\begin{Thm}
The embedding  $N-D \rightarrow N_{\gamma}$ induces an isomorphism
$$S(N-D)\cong S(N_{\gamma}).$$
\end{Thm}
\begin{proof}
We proceed by induction on $n,$ the number of the $2$-handles to be added to $({{\cal H}}_{m})_{\gamma}$ to obtain $N_{\gamma}.$ If $n=0,$ we are done by Theorem 5. If $n \ge 1,$
let $N'$ be the $3$-manifold obtained from $({{\cal H}}_{m})_{\gamma}$ by adding $(n-1)$ of those $2$-handles added to $({{\cal H}}_{m})_{\gamma}$. Suppose the result is true for $N'$, i.e.
$$S(N'-D)\cong S(N'_{\gamma}).$$

Suppose that the $n$th $2$-handle is added along a curve $\gamma^{*}$ in the boundary of $({{\cal H}}_{m}),$  where $\gamma^{*}$ is disjoint from $\gamma$ and the  curves where the other $(n-1)$ $2$-handles are attached. Let $A'$ and $B'$ be two points on $\gamma^{*}.$ By the proof of the Epimorphism Lemma 2.2,
$$S(N'-D, A', B') \twoheadrightarrow S(N'_{\gamma}, A', B').$$
Using \cite[section 2]{ZG99}, we have the following commutative diagram with exact rows.
$$
\begin{CD}
S(N'-D, A', B')@>>> 
S(N'-D) @>>>  S({(N'-D)}_{\gamma^{*}}) @>>> 0
\\
@V\text{onto} VV   @V\cong VV   @V VV
\\
S(N'_{\gamma}, A', B')@>>> 
S(N'_{\gamma}) @>>>  S({(N'_{\gamma})}_{\gamma^{*}}) @>>> 0
\end{CD}
$$
\noindent The vertical map on the right is an isomorphism by the five-lemma. $N$ is  obtained from $N'$ by adding the $n$th $2$-handle along $\gamma^{*}$.
Thus $(N'-D)_{\gamma^{*}} =N-D$ and $(N'_{\gamma})_{\gamma^{*}}=N_{\gamma}$.
\end{proof}

\begin{Cor} Let $B_1$ and $B_2$ denote the 3-balls we remove from $M_1$ and $M_2$ while forming $M_1\# M_2.$ The embedding  $(M_1 - B_1)\sqcup (M_2 - B_2) \rightarrow M_1\# M_2$ induces an isomorphism
$$S(M_1)\otimes S(M_2)\cong S(M_1\# M_2).$$
\end{Cor}
\begin{proof}
Since $S(M-D) \cong S(M_1)\otimes S(M_2)$.
\end{proof}

The above corollary holds whether or not $M_1$ or $M_2$ have boundary.

\section {An example in $S^1\times S^2\# S^1\times S^2$}

In \cite{ZG99}, we showed that $S(S^1\times S^2)$ is a free ${\cal R}[x, x^{-1}]$-module generated by the empty link. It follows that $S(S^1\times S^2\# S^1\times S^2)$ is also a free module generated by the empty link.
Let $K$ be a knot in $S^1\times S^2\# S^1\times S^2$ pictured by the following diagram:\vglue-10mm
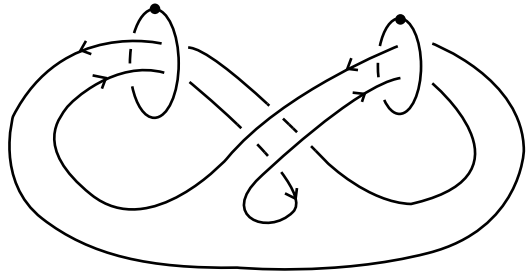
$$\raisebox{0mm}{\epsfxsize.0in\epsffile{exs1s2.eps}}$$
Here the two circles  with a dot are a framed link description of $S^1\times S^2\# S^1\times S^2$. 
Note this same knot was studied with respect to the Kauffman Bracket skein modules in 
\cite{pG98}.

In $S(S^2\times S^1 \# D^3, 4pts)$, isotopy  yields,
$$\raisebox{-12mm}{\epsfxsize.0in
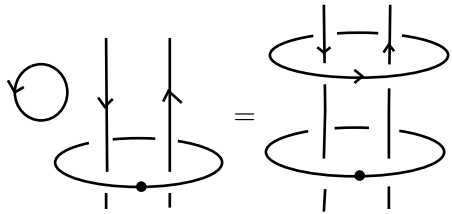
\epsffile{ex1.eps}}= \raisebox{-12mm}{\epsfxsize.0in\epsffile{ex2.eps}}$$
Using the Homflypt skein relations in $S(D^2\times I, 4pts)$,
$$\raisebox{-7mm}{\epsfxsize.0in\epsffile{ex21.eps}}=\Big(\dfrac{v^{-1}-v}{s-s^{-1}}-(v-v^{-1})(s-s^{-1})\Big)\raisebox{-7mm}{\epsfxsize.0in\epsffile{ex3.eps}}-(s-s^{-1})^2\raisebox{-7mm}{\epsfxsize.0in\epsffile{ex4.eps}}$$
Therefore, in $S(S^2\times S^1 \# D^3, 4pts)$, we have: 
$$(v-v^{-1})\Big)\raisebox{-7mm}{\epsfxsize.0in\epsffile{ex5.eps}}\equiv -(s-s^{-1})\raisebox{-11mm}{\epsfxsize.0in\epsffile{ex6.eps}}$$
Thus\vglue-10mm
$$(v-v^{-1})^2 K= (s-s^{-1})^2\raisebox{-20mm}{\epsfxsize.0in\epsffile{exs1s21.eps}}$$
$=(s-s^{-1})^2 \frac{v^{-1}-v}{s-s^{-1}}\phi. $
i.e. $K=\dfrac{s-s^{-1}}{v^{-1}-v}\phi$ in $S(S^1\times S^2\# S^1\times S^2)$.

\section{Proof of Lemma 2.1}

Note $y_{\lambda}=F_{\lambda}G_{\lambda}$. We start with the following:
$$\raisebox{-18mm}{\epsfxsize.0in\epsffile{albesi01.eps}}=\raisebox{-18mm}{\epsfxsize.0in\epsffile{albesi0.eps}}=$$
\begin{equation}
x^2\raisebox{-18mm}{\epsfxsize.0in\epsffile{albesia1.eps}}+x(s-s^{-1})\raisebox{-18mm}{\epsfxsize.0in\epsffile{albesia2.eps}}\tag{I*}
\end{equation}
$$=x^2\raisebox{-18mm}{\epsfxsize.0in\epsffile{albesia1.eps}}+x(s-s^{-1})(xv^{-1})\raisebox{-18mm}{\epsfxsize.0in\epsffile{albea2.eps}}$$
We pulled out the string corresponding to the last cell in the last row of $\lambda.$ Therefore in the above diagram, a $1$ by the side of the string indicates the string related to the last cell in the last row of $\lambda$.  Applying the Homflypt skein relation to the last diagram: 
$$\raisebox{-18mm}{\epsfxsize.0in
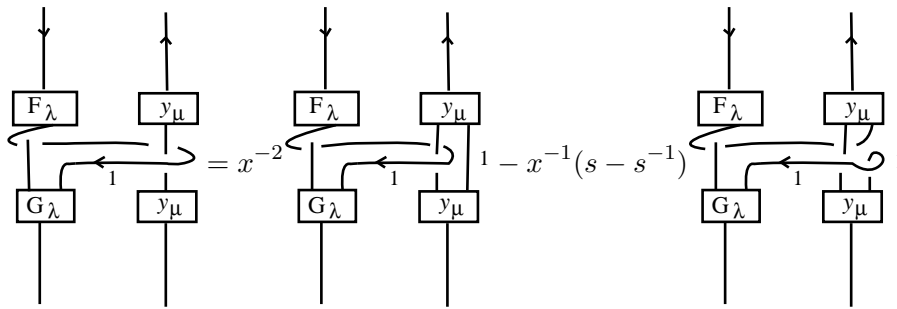
\epsffile{albea2.eps}}=x^{-2}\raisebox{-18mm}{\epsfxsize.0in\epsffile{albea21.eps}}-x^{-1}(s-s^{-1})\raisebox{-18mm}{\epsfxsize.0in\epsffile{albea22.eps}}$$
We pulled out the string corresponding to the last cell in the last row of $\mu.$ Continuing to pull out strings which correspond to cells of $\mu,$ working to the left through columns and upward through the rows of $\mu,$ we obtain:

$$\raisebox{-18mm}{\epsfxsize.0in
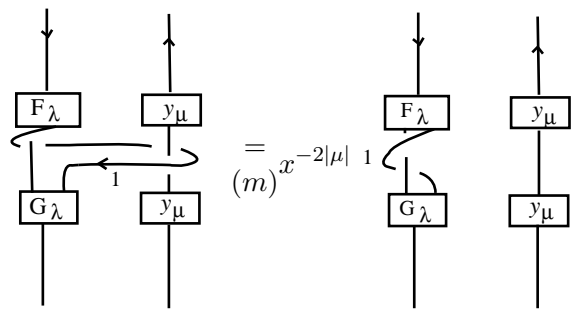
\epsffile{albea2.eps}}
\quad
\begin{matrix} = \\{(m)} \end{matrix}
x^{-2|\mu|}\ \raisebox{-18mm}{\epsfxsize.0in\epsffile{albea1.eps}}$$
where the string corresponding to the last cell in the last row of $\lambda$ encircles the remaining $|\lambda|-1$ strings as shown.

In this way Equation (I*) becomes:

$$\raisebox{-18mm}{\epsfxsize.0in
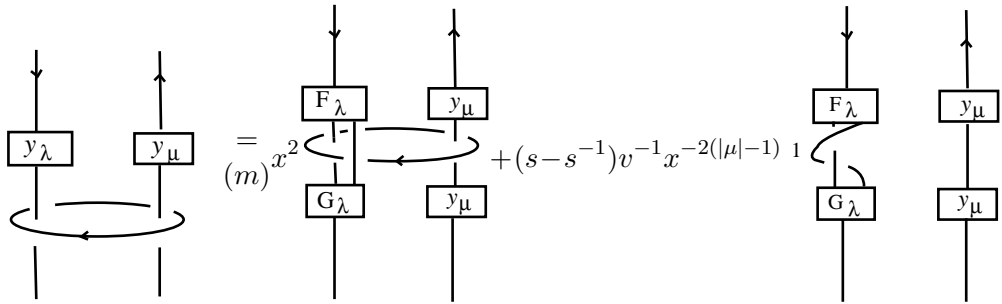
\epsffile{albesi01.eps}}
\quad
\begin{matrix} = \\{(m)} \end{matrix}
x^2\raisebox{-18mm}{\epsfxsize.0in\epsffile{albesia1.eps}}
+(s-s^{-1})v^{-1}x^{-2(|\mu|-1)}
\ \raisebox{-18mm}{\epsfxsize.0in\epsffile{albea1.eps}}$$

We continue in this way, pulling the encircling component successively through the vertical strings corresponding to cells of $\lambda$, working to the left through columns and upward through the rows of $\lambda.$  We obtain:
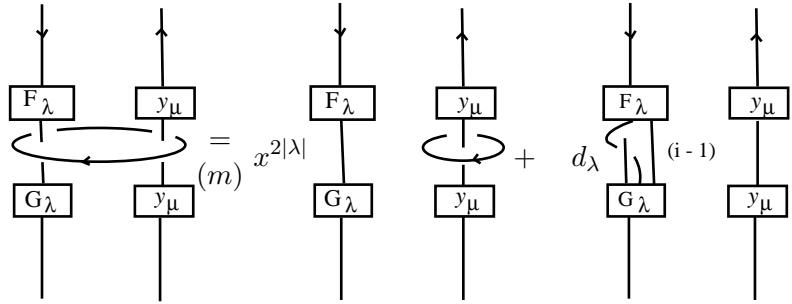
\vglue -5mm
$$
 \raisebox{-18mm}{\epsfxsize.0in\epsffile{albesi0.eps}}
\begin{matrix} = \\{(m)} \end{matrix} \
x^{2|\lambda|}\raisebox{-18mm}{\epsfxsize.0in\epsffile{albemu0.eps}}
+
\quad d_\lambda \raisebox{-18mm}{\epsfxsize.0in\epsffile{albeai.eps}}
$$

where $d_\lambda$ denotes $v^{-1}(s-s^{-1})\sum_{i=1}^{|\lambda|}x^{-2(|\mu|-i)}.$
In the last diagram, the  $i-1$ vertical strings are related to the last $i-1$ cells of $\lambda$ by the index order, and the $i$th string encircles the remaining $|\lambda|-i$ strings. Lemma 2.1 follows from the following lemma and Lemma 6.2 (a) below.

\begin{Lem} Let $\mu$ be a Young diagram of size $n$,
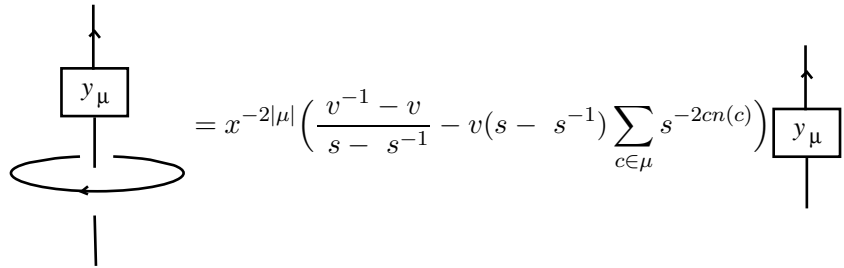
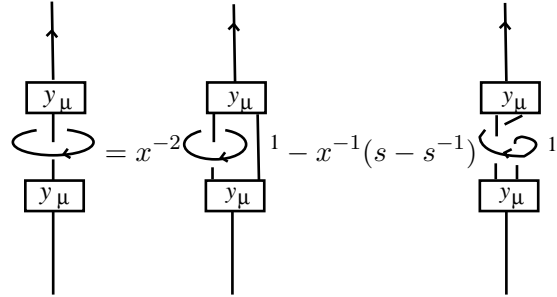
$$\raisebox{-18mm}{\epsfxsize.0in\epsffile{albeym.eps}}= x^{-2|\mu|}\Big(\dfrac{v^{-1}-v} {\ s - \ s^{-1}}-v(s - \ s^{-1})\sum_{c \in \mu}s^{-2cn(c)}\Big)\raisebox{-10mm}{\epsfxsize.0in\epsffile{albeym1.eps}}$$
\end{Lem}
\begin{proof}
We consider
$$\raisebox{-18mm}{\epsfxsize.0in\epsffile{albmu0.eps}}=x^{-2}\raisebox{-18mm}{\epsfxsize.0in\epsffile{albemu01.eps}}-x^{-1}(s-s^{-1})\raisebox{-18mm}{\epsfxsize.0in\epsffile{albem011.eps}}$$
Here we start with the string corresponding to the last cell in the last row of $\mu$, we pull the encircling component successively through the vertical strings, working to the left through columns and upward through the rows. Repeating the above process, for $i\geq 2$:
$$=x^{-2i}\raisebox{-18mm}{\epsfxsize.0in
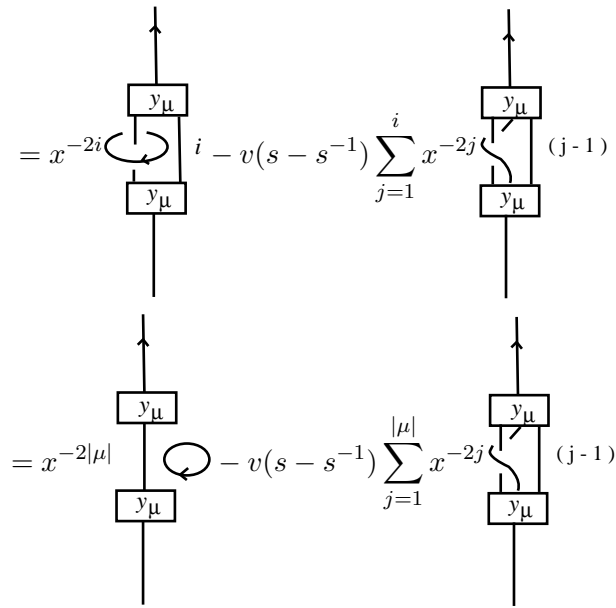
\epsffile{albemu0i.eps}}-v(s-s^{-1})\sum_{j=1}^{i}x^{-2j}\raisebox{-18mm}{\epsfxsize.0in\epsffile{albem0i1.eps}}$$
$$=x^{-2|\mu|}\raisebox{-18mm}{\epsfxsize.0in\epsffile{albe00.eps}}-v(s-s^{-1})\sum_{j=1}^{|\mu|}x^{-2j}\raisebox{-18mm}{\epsfxsize.0in\epsffile{albem0i1.eps}}$$
The result follows from Lemma 6.2 (b) below.
\end{proof}

\begin{Lem} Let $\lambda$ be a Young diagram and $(h, l)$ be the index of the cell after which $(i-1)$ cells of $\lambda$ follow. 

(a)\vglue -10mm
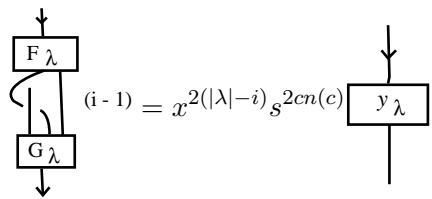
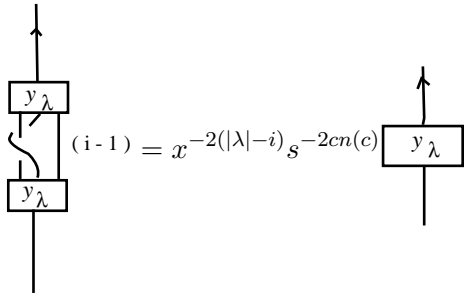
$$\raisebox{-11mm}{\epsfxsize.0in\epsffile{f1.ai}}=x^{2(|\lambda|-i)}s^{2cn(c)}\raisebox{-9mm}{\epsfxsize.0in\epsffile{fig22.eps}} $$
(b)\vglue -10mm
$$\raisebox{-18mm}{\epsfxsize.0in\epsffile{f2.ai}}=x^{-2(|\lambda|-i)}s^{-2cn(c)}\raisebox{-9mm}{\epsfxsize.0in\epsffile{fig21.eps}}$$
\end{Lem}
\noindent \remark The techniques used in this proof are similar to the proof of the framing factor in section 5 of \cite{AM98} by H. Morton and A. Aiston.
\begin{proof}
(a) We will borrow the notation of H. Morton and A. Aiston and use a schematic dot diagram to represent the element in the Hecke category $H_{\square_{\lambda}}$, which is between $F_{\lambda}$ and $G_{\lambda}$ as shown on the left-hand side. 

Recall that $y_{\lambda}=F_{\lambda}G_{\lambda}$. Now in the diagram of the left-hand side of (a), introduce a schematic picture $T$ as follows:
$$T=\raisebox{-15mm}{\epsfxsize.0in
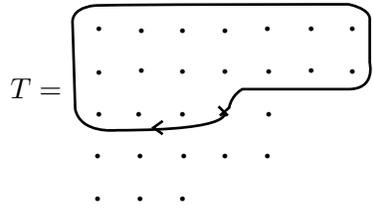
\epsffile{amdot.eps}}$$
This indicates that the last $i-1$ strings were pulled out, the $i$th string marked by $\times$ starts and finishes at $(h, l)$. The arrow on the $i$th string shows the string orientation when we look at it from above. The $i$th string encircles the remaining $|\lambda|-i$ strings in the clockwise direction. Here all strings shown by single dots are going vertical. 
The left-hand side of (a) can be expressed as $F_{\lambda}TG_{\lambda}$. We will be working on $F_{\lambda}TG_{\lambda}$. Using the Homflypt skein relations and the inseparability in Lemma 16 of \cite{AM98}, we have, $$F_{\lambda}TG_{\lambda}=x^{2(|\lambda|-(i-1)-hl)}F_{\lambda}SG_{\lambda}$$
Where $S$ is given by:
$$S=\raisebox{-10mm}{\epsfxsize.0in
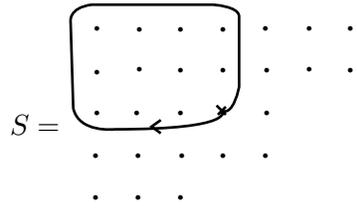
\epsffile{asdot.eps}}$$
Since $S=S_{1}T_{1}S_{2}$, where:
$$S_1=\raisebox{-15mm}{\epsfxsize.0in
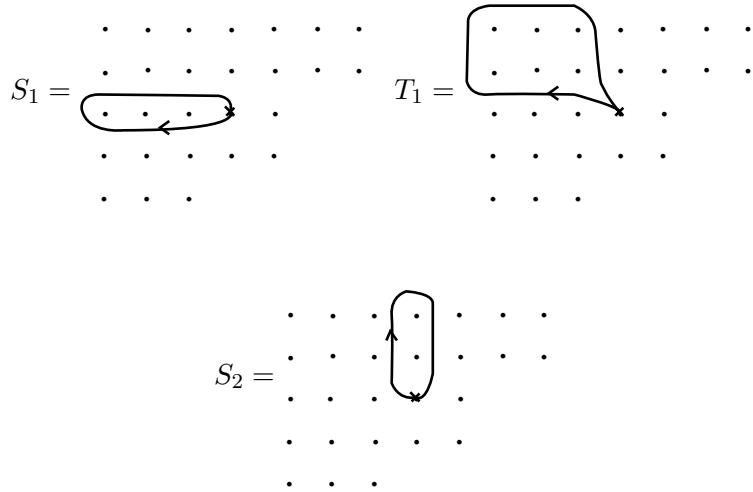
\epsffile{as1dot.eps}}\quad T_1=\raisebox{-15mm}{\epsfxsize.0in\epsffile{at1dot.eps}}$$
$$S_2=\raisebox{-15mm}{\epsfxsize.0in\epsffile{as2dot.eps}}$$
First we have $F_{\lambda}S_1=(xs)^{2(l-1)}F_{\lambda}$ by the property $\sigma_{i}f_{m}=xsf_{m}$ ; secondly, $S_{2}G_{\lambda}=(-xs^{-1})^{2(h-1)}G_{\lambda}$ by the property $g_{m}\sigma_{i}=-xs^{-1}g_{m}$,  \cite[Lemma 8]{AM98}) .
It follows that   $F_{\lambda}SG_{\lambda}=x^{2(h+l-2)}s^{2(l-h)}F_{\lambda}T_{1}G_{\lambda}$. By a similar argument as in the proof of Theorem 17 in \cite{AM98}, $F_{\lambda}T_{1}G_{\lambda}=x^{2(l-1)(h-1)}F_{\lambda}G_{\lambda}$. Thus $F_{\lambda}TG_{\lambda}=x^{2(|\lambda|-i)}s^{2(l-h)}F_{\lambda}G_{\lambda}=x^{2(|\lambda|-i)}s^{2cn(c)}y_{\lambda}$, where c is the cell indexed by $(h,l)$.

(b) We prove the result with all string orientations reversed. As string reversal  defines a skein module isomorphism, this suffices. As $y_{\lambda}=F_{\lambda}G_{\lambda}$, we can use the following schematic picture to denote the left-hand side of (b) as $F_{\lambda}G_{\lambda}T^{-1}F_{\lambda}G_{\lambda}$, where
$$T^{-1}=\raisebox{-15mm}{\epsfxsize.0in
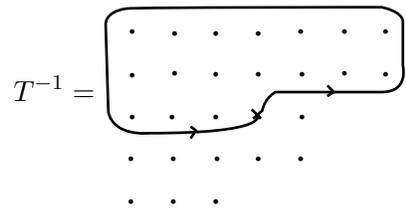
\epsffile{aldot.eps}},$$
the $i$th string is indexed by $(h, l)$ and circles the remaining strings in the clockwise direction. Again, we have
$$G_{\lambda}T^{-1}F_{\lambda}=x^{-2(|\lambda|-(i-1)-hl)}G_{\lambda}S^{-1}F_{\lambda}.$$
Where $S^{-1}$ is given by:
$$S^{-1}=\raisebox{-10mm}{\epsfxsize.0in\epsffile{as11dot.eps}}$$
Since $S^{-1}=S_{2}^{-1}T_{1}^{-1}S_{1}^{-1}$, where:
$$S_{2}^{-1}=\raisebox{-15mm}{\epsfxsize.0in\epsffile{as21dot.eps}},\quad T_{1}^{-1}=\raisebox{-15mm}{\epsfxsize.0in\epsffile{att1dot.eps}}$$
$$S_{1}^{-1}=\raisebox{-15mm}{\epsfxsize.0in\epsffile{as22dot.eps}}.$$
We have $G_{\lambda}S_{2}^{-1}=(-x^{-1}s)^{2(h-1)}G_{\lambda}$ and $S_{1}^{-1}F_{\lambda}=(x^{-1}s^{-1})^{2(l-1)}F_{\lambda}$ by the properties $\sigma_{i}f_{m}=xsf_{m}$ and $g_{m}\sigma_{i}=-xs^{-1}g_{m}$.

We have $$G_{\lambda}S^{-1}F_{\lambda}=x^{-2(h+l-2)}s^{2(h-l)}G_{\lambda}T_{1}^{-1}F_{\lambda}=x^{-2(h+l-2)}s^{2(h-l)}x^{-2(h-1)(l-1)}G_{\lambda}F_{\lambda}.$$ 
It follows that  
$G_{\lambda}T^{-1}F_{\lambda}=x^{-2(|\lambda|-i)}s^{-2cn(c)}G_{\lambda}F_{\lambda}$. By the idempotent property, $F_{\lambda}G_{\lambda}T^{-1}F_{\lambda}G_{\lambda}=x^{-2(|\lambda|-i)}s^{-2cn(c)}F_{\lambda}G_{\lambda}$. The result follows.
\end{proof}
\section {Discussion of the proofs of Theorems 2 \& 3}

The proof of Theorem 2 is  basically the same as the proof of Theorem 1. 
However as noted in the introduction we do not yet know that the relative Homflypt skein of a handlebody is  free. So we must work over $\cal {F}.$

For the proof of Theorem 3, we note that  every relative  link in $(M,X,Y)$ is isotopic to a link which intersects a tubular neighborhood of $\Sigma$ with
 $m$ straight strands going in one direction and $m+r$
straight strands going the other direction.
We will write such elements as linear combinations of $\alpha_t\beta_{\tau}\otimes {\alpha_{\sigma}\beta_{s}}^{*}$, and $t$ and $\tau$ are standard tableaux of a Young diagram $\lambda$, $\sigma$ and $s$ are standard tableaux of a Young diagram $\mu$ with
 $|\lambda|=m,$ and $|\mu|=m+r$. As $x^{2r}-1- c_{\lambda,\mu}$ is invertible over $k_r,$ we have that $\alpha_t\beta_{\tau}\otimes { \alpha_{\sigma}\beta_{s} }^{*}
\in {\cal L}(|\lambda|+|\mu|).$  We may repeat this argument until the  class of our original relative link is represented by a linear combination of links  each of which intersects $\Sigma$ less than $r$ times. This must be the empty linear combination.

{}
\Addresses


\begin{thebibliography}{10}


\bibitem{AM98}
A.~K.~Aiston and H.~R.~Morton, {\em Idempotents of Hecke algebras of type A}, J. of Knot Theory and Ram. 7 No 4 (1998), 463-487.

\bibitem{cB98}
C.~Blanchet, {\em Hecke algebras, modular categories and $3$-manifolds quantum invariants}, Topology (39) 1 (2000), 193-223.

\bibitem{pG98} P.~Gilmer, {\em A TQFT for wormhole cobordisms over the field of rational functions} in Knot theory (Warsaw, 1995), 119--127, Banach Center Publ., 42, Polish Acad. Sci., Warsaw, 1998.

\bibitem{ZG99}
P.~Gilmer and J.~K.~Zhong, {\em The Homflypt skein module of $S^{1}\times S^{2}$}, Math Zeit., vol 237, pp 769-814 (2001).

\bibitem{J}
N. ~Jacobson, {\em Basic Algebra II}, second edition,
 W.H. Freeman (1989).

\bibitem{M93}
H.~Morton, {\em Invariants of links and 3-manifolds from skein theory and from quantum groups},
in {\em Topics in knot theory.} N.A.T.O. A.S.I. series C 399, eds. M. ~Bozhuyuk, Kluwer (1993) 479-497.

\bibitem{MT90}
H.~Morton and P.~Traczyk, {\em Knots and algebras}, Contribuciones Matematicas en homaje al professor D.~Antonio Plans Sanz de Bremond, E.~Martin-Peinador and A.~Rodez editors, University of Saragoza (1990), 201-220.

\bibitem{JP90} 
J.~Przytycki, {\em Skein module of links in a handlebody},
Topology '90 (Columbus, OH, 1990), 315--342, 
Ohio State Univ. Math. Res. Inst. Publ., 1, 
de Gruyter, Berlin, 1992.


\bibitem{JP} 
 J.~Przytycki, {\em Kauffman bracket skein module of a connected sum of $3$-manifolds},  Manuscripta Math. 101 (2000), no. 2, 199--207.



\bibitem{vT90}
V.~G.~Turaev, {\em Operator invariants of tangles, and R-matrices}, Math. USSR Izv. Vol. 35 No. 2 (1990), 411-443.


 

\end{thebibliography}
\end{document}